\documentclass[aop]{imsart}
\usepackage[ansinew]{inputenc}
\usepackage{amsthm,amsmath,amssymb}
\usepackage{enumerate}
\RequirePackage[colorlinks,citecolor=blue,urlcolor=blue]{hyperref}

\newcommand\NoBlackBoxes{\global\overfullrule0pt}
\NoBlackBoxes
\theoremstyle{plain} 

\def\4{\kern1pt}

\def\6{\vphantom0}

\def\8{\kern-10pt}
\def\7#1{_{(#1)}}

\makeatletter
\let\serieslogo@\relax
\let\@setcopyright\relax

\def\speciallabelmark#1{\def\@currentlabel{#1}}
\makeatother

\begin{document}

\def\ffrac#1#2{\raise.5pt\hbox{\small$\4\displaystyle\frac{\,#1\,}{\,#2\,}\4$}}
\def\ovln#1{\,{\overline{\!#1}}}
\def\ve{\varepsilon}
\def\kar{\beta_r}
\def\theequation{\thesection.\arabic{equation}}
\def\E{{\mathbb E}}
\def\R{{\mathbb R}}
\def\C{{\bf C}}
\def\P{{\mathbb P}}
\def\H{{\rm H}}
\def\Im{{\rm Im}}
\def\Tr{{\rm Tr}}

\def\k{{\kappa}}
\def\M{{\cal M}}
\def\Var{{\rm Var}}
\def\Ent{{\rm Ent}}
\def\O{{\rm Osc}_\mu}

\def\ep{\varepsilon}
\def\phi{\varphi}
\def\F{{\cal F}}
\def\L{{\cal L}}

\def\s{{\mathfrak s}}

\def\be{\begin{equation}}
\def\en{\end{equation}}
\def\bee{\begin{eqnarray*}}
\def\ene{\end{eqnarray*}}

\begin{frontmatter}
	
\title{NORMAL APPROXIMATION FOR WEIGHTED SUMS\\
		UNDER A SECOND ORDER CORRELATION CONDITION}
\runtitle{CLT for Weighted Sums}

\begin{aug}
\author{\fnms{S. G.} \snm{Bobkov}\thanksref{t1,m1}\ead[label=e1]{bobkov@math.umn.edu}},
\address{School of Mathematics\\ University of Minnesota \\ 
    Vincent Hall 228\\
		206 Church St. S.E.\\
		Minneapolis, MN 55455\\ USA\\
\printead{e1}}
\author{\fnms{G. P.} \snm{Chistyakov}\thanksref{t1,m2}\ead[label=e3]
{chistyak@math.uni-bielefeld.de}}
\address{Department of Mathematics\\
		University of Bielefeld\\
		Postbox 100131\\
		33501 Bielefeld\\
		Germany\\\printead{e2}\\
\phantom{E-mail:\ }
\printead*{e3}}
	\and
\author{\fnms{F.} \snm{G\"otze}
\thanksref{t1,m2}\ead[label=e2]{goetze@math.uni-bielefeld.de}}
	
\thankstext{t1}{Supported by NSF grant DMS-1855575, and CRC 1283 at 
Bielefeld University.}

\runauthor{S. Bobkov, C.  Chistyakov and F. G\"otze}
\affiliation{University of Minnesota\thanksmark{m1} and 
Bielefeld University\thanksmark{m2}}
\end{aug}

	
	
	
	
	

%
%
%



\begin{abstract}
\hskip-5.5mm
Under correlation-type conditions, we derive an upper bound of order 
$(\log n)/n$ for the average Kolmogorov distance between the distributions 
of weighted sums of dependent summands and the normal law. The result
is based on improved concentration inequalities on high-dimensional Euclidean 
spheres. Applications are illustrated on the example of log-concave probability 
measures.
\end{abstract}



\end{frontmatter}


\section{{\bf Introduction}}
\setcounter{equation}{0}

Let $X = (X_1,\dots,X_n)$ be an isotropic random vector in $\R^n$ 
($n \geq 2$), that is, with uncorrelated components having mean zero 
and variance one. We consider the distribution functions 
$F_\theta(x) = \P\{S_\theta \leq x\}$ of the weighted sums
$$
S_\theta = \theta_1 X_1 + \dots + \theta_n X_n, \qquad 
\theta = (\theta_1,\dots,\theta_n), \ \ \theta_1^2 + \dots + \theta_n^2=1,
$$
with coefficients taken from the unit sphere $S^{n-1}$ in $\R^n$.
Thus, $\E S_\theta = 0$ and $\Var(S_\theta) = 1$ for all
$\theta \in S^{n-1}$.

The central limit problem is to determine natural conditions on $X$ and 
$\theta$ which ensure that the random variables $S_\theta$ are nearly 
standard normal. In this case, one would also like to explore the rate 
of normal approximation in the Kolmogorov distance
$$
\rho(F_\theta,\Phi) = \sup_x |F_\theta(x) - \Phi(x)|,
$$
where 
$$
\Phi(x) = \frac{1}{\sqrt{2\pi}} \int_{-\infty}^x e^{-y^2/2}\,dy, 
\qquad x \in \R,
$$
is the standard normal distribution function. Let us briefly recall 
several well-known results in the case of independent components $X_k$.
Here, one of general variants of the central limit theorem asserts that 
$\rho(F_\theta,\Phi)$ will be small, as long as $X_k$ are identically 
distributed (the i.i.d. case), while $\max_k |\theta_k|$ is small.
Moreover, under the 3rd moment condition this property may be
quantified by virtue of the Berry-Esseen bound
\be
\rho(F_\theta,\Phi) \leq c \sum_{k=1}^n |\theta_k|^3\,\E\,|X_k|^3.
\en
Here and below, we denote by $c$, or by $c_j$ with an integer index $j$
absolute positive constants which may vary from place to place.
The inequality (1.1) extends to the non-i.i.d. case as well (\cite{P1}, \cite{P2}). 

It easy to see that the sum in (1.1) is greater than or equal to $1/\sqrt{n}$ 
for all $\theta$, and that (1.1) leads to this standard 
$\frac{1}{\sqrt{n}}$-rate in the i.i.d. case, once the coefficients $\theta_k$ 
are equal to each other. For general distributions of $X_k$ this 
standard rate cannot be improved by assuming stronger moment-type conditions. 
Nevertheless, one may look at the problem from an ensemble point of view in 
$\theta$ asking whether or not $\rho(F_\theta,\Phi)$ will be essentially 
smaller than $1/\sqrt{n}$ for most of $\theta$ on the sphere measured with 
the uniform probability measure $\s_{n-1}$ on $S^{n-1}$. A striking result 
in this direction was obtained by Klartag and Sodin \cite{K-S}, showing 
in particular that
\be
\E_\theta \rho(F_\theta,\Phi) \leq \frac{c}{n}\,\bar \beta_4, \qquad
\bar \beta_4 = \frac{1}{n} \sum_{k=1}^n \E X_k^4,
\en
where we use $\E_\theta$ to denote the average over the measure $\s_{n-1}$.
Large deviation bounds for the set on the sphere where $\rho(F_\theta,\Phi)$ 
exceeds a multiple of $\frac{1}{n}\,\bar \beta_4$ are derived in \cite{K-S} as well. 
Thus, when $\bar \beta_4$ is bounded like in the i.i.d. case, 
the distances $\rho(F_\theta,\Phi)$ turn out to be typically 
of order $1/n$ in contrast to the classical case of equal coefficients.

The aim of these notes is to extend this interesting phenomenon under 
a suitable correlation-type condition (and thus for some class of dependent 
$X_k$) to isotropic random vectors with a similar $\frac{1}{n}$-rate modulo 
a logarithmic factor. The scheme of the weighted sums under dependence has 
already a long history, going back to the seminal work of Sudakov \cite{Su}. 
We will give a short overview of this line of research in Section 10
(partly in Section 7), and now turn to the main result.

We will say that the random vector $X$ satisfies a second order correlation 
condition with constant $\Lambda$, if for any collection $a_{ij} \in \R$,
\be
\Var\bigg(\sum_{i,j=1}^n a_{ij} X_i X_j\bigg) \leq 
\Lambda \sum_{i,j=1}^n a_{ij}^2.
\en
An optimal value $\Lambda = \Lambda(X)$ is finite as long as $|X|$ has a finite 
$4$-th moment, and then it represents the maximal eigenvalue of the covariance 
matrix associated with the $n^2$-dimensional random vector 
$\big(X_i X_j - \E X_i X_j\big)_{i,j=1}^n.$

\vskip5mm
{\bf Theorem 1.1.} {\sl Let $X$ be an isotropic random vector in $\R^n$
with a symmetric distribution and a finite constant $\Lambda = \Lambda(X)$. 
Then 
\be
\E_\theta\, \rho(F_\theta,\Phi) \leq \frac{c\,\log n}{n}\, \Lambda.
\en
}

The characteristic $\Lambda$ may be bounded, for example, via the relation
$\Lambda \leq 4/\lambda_1$ in terms
of a positive spectral gap, that is in terms of the optimal value 
$\lambda_1 = \lambda_1(X)$ in the Poincar\'e-type inequality
\be
\lambda_1\, \Var(u(X)) \, \leq \, \E\, |\nabla u(X)|^2
\en
(with $\lambda_1>0$), where $u$ is an arbitrary smooth function $u$ on $\R^n$
(cf. Proposition 3.4 below). In one important particular case, the well-known 
Kannan-Lov\'asz-Simonovits conjecture asserts that 
$\lambda_1$ is bounded away from zero for the whole class of isotropic 
log-concave probability distributions on the Euclidean space $\R^n$ of 
any dimension (for short, K-L-S). Conditional on K-L-S, Theorem 1.1 would 
hence guarantee the $\frac{\log n}{n}$-rate.

\vskip5mm
{\bf Corollary 1.2.} {\sl Let $X$ be an isotropic random vector in $\R^n$
with a symmetric log-concave distribution. Assuming the K-L-S hypothesis,
we have
\be
\E_\theta\, \rho(F_\theta,\Phi) \leq \frac{c\,\log n}{n}.
\en
}

In fact, modulo a logarithmic factor, the conclusion may be reversed in the 
sense that (1.6) implies $1/\lambda_1 \leq c\,(\log n)^7$, cf. Section 8.

An unconditional statement in the isotropic log-concave case with 
a standard rate of normal approximation can be obtained by combining
the results of \cite{A-B-P} and \cite{B1} on the concentration of $F_\theta(x)$
around the average distribution function $F(x) = \E_\theta F(x)$
with respect to the variable $\theta$ with a recent bound
in the thin-shell problem due to Lee and Vempala \cite{L-V} on the
concentration of the Euclidean length $|X|$ about its average value $\E\,|X|$
(which is in essense equivalent to the closeness of $F$ to the standard normal
distribution function $\Phi$). More details are given in Section 7; one then gets
\be
\E_\theta\, \rho(F_\theta,\Phi) \leq \frac{c}{\sqrt{n}}\sqrt{\log n}.
\en

As for the general (not necessarily log-concave) case, the functional $\Lambda(X)$
turns out to be responsible for both, formally different concentration problems.
The proof of Theorem 1.1 is based on results for spherical concentration, which 
have been recently developed in \cite{B-C-G1}. They provide improved rates of 
concentration for smooth functions $u$ on the sphere based on the additional 
information about the Hessian of $u$. This naturally leads 
to the definition of $\Lambda(X)$ as introduced above. The ``2nd-order'' 
concentration inequalities on $S^{n-1}$ may also be used to derive 
large deviation bounds for $\rho(F_\theta,\Phi)$ considered as random variables 
on the probability space $(S^{n-1},\s_{n-1})$.
Moreover, one may remove the symmetry assumption as well, by adding to the
right-hand side of (1.4) an additional term responsible for 3rd order
correlations between $X_k$. We refrain from including these somewhat more
technical results here and refer the interested reader to \cite{B-C-G4} for 
a full account. 

As we will see, there exist several natural classes of 
probability distributions for which a bound on the parameter $\Lambda$ 
can be obtained. Some of them are considered in Section 3, after a brief
discussion of general properties of $\Lambda$ and related functionals
in Section 2. Some results about the second order concentration on the sphere
are described in Sections 4, which we apply in Section 5 to explore the
concentration of characteristic functions of $S_\theta$ with respect
to the variable $\theta$. In Section 6, relying upon a general Berry-Esseen-type 
inequality, we finalize the proof of Theorem 1.1.
The relationship of Theorem 1.1 with the K-L-S conjecture and a closely
related thin-shell problem in the log-concave case are discussed separately
in Sections 7-8.


\vskip5mm
\section{{\bf Second Order Correlation Condition and Related Functionals}}
\setcounter{equation}{0}

As usual, the Euclidean space $\R^n$ is endowed with the canonical norm 
$|\,\cdot\,|$ and the inner product $\left<\cdot,\cdot\right>$. We start with 
preliminary remarks on the second order correlation condition and related 
functionals.

Let $X = (X_1,\dots,X_n)$ be a random vector in $\R^n$. With the 
Hilbert-Schmidt norm of a matrix $A = (a_{ij})_{i,j=1}^n$ given by 
$\|A\|_{\rm HS} = (\sum a_{ij}^2)^{1/2}$, the definition (1.3) becomes
$$
\Var\big(\left<AX,X\right>\big) \leq \Lambda \|A\|_{\rm HS}^2,
$$
where we may restrict ourselves to symmetric 
matrices $A$ only. This description shows that the functional $\Lambda(X)$ 
is invariant under linear orthogonal transformations of the space $\R^n$ 
(just as the Hilbert-Schmidt norm).

Related moment and variance-type functionals are
\bee
M_p 
 & = &
M_p(X) \ = \
\sup_{\theta \in S^{n-1}}\, (\E\, |S_\theta|^p)^{1/p} \ \ \ (p \geq 1), \\
\sigma_4^2 
 & = &
\sigma_4^2(X)  \ = \ \frac{1}{n}\,\Var(|X|^2).
\ene
We are mostly interested in the moments $M_p$ with $p=2$ and $p=4$.
For example, $M_2 = 1$ in the isotropic case, and $\sigma_4 = 0$, if
$|X|$ is constant a.s. These functionals can be controlled in terms of 
$\Lambda$, as the following statement shows.

\vskip2mm
{\bf Proposition 2.1.} {\sl We have
$$
a) \ M_4^4 \leq M_2^4 + \Lambda; \qquad  b) \ \sigma_4^2 \leq \Lambda.
$$
}

\vskip2mm
{\bf Proof.} Choosing in (1.3) $a_{ij} = \theta_i \theta_j$, 
$\theta \in S^{n-1}$, we get $\Var(S_\theta^2) \leq \Lambda$.
Since $\E S_\theta^2 \leq M_2^2$, it follows that 
$\E S_\theta^4 \leq M_2^4 + \Lambda$, that is, $a)$.
Putting $a_{ij} = \delta_{ij}$, we also obtain $b)$.
\qed

\vskip5mm
In turn, the $M_p$-moments may be related to the moments of $|X|$. 
It is easy to see that
$$
(\E\, |X|^p)^{1/p} \leq M_p \sqrt{n}, \qquad p \geq 2,
$$
while in the isotropic case, there is an opposite inequality
$(\E\, |X|^p)^{1/p} \geq \sqrt{n}$.

The functionals $\sigma_4^2$, $M_4$, and $\Lambda$ are
useful for the estimation of ``small" ball probabilities. 
For example, if $\E\,|X|^2 = n$, using an independent copy $Y$ of $X$,
we have
$$
\P\Big\{|X - Y|^2 \leq \frac{1}{4}\,n\Big\} \, \leq \, 
\frac{A}{n^2}, \qquad A = 256\,(M_4^8 + \sigma_4^4).
$$
This bound was applied in the proof of Lemma 5.1 below (for details 
we refer to \cite{B-C-G3}). Here, by Proposition 2.1 $a)$ in the isotropic case, 
$A \leq c\Lambda^2$, which is also due to 
the fact that the functional $\Lambda(X)$ is bounded away 
from zero for $n \geq 2$ (in contrast to $\sigma_4$).

\vskip5mm
{\bf Proposition 2.2.} {\sl If $X$ is isotropic, then 
$\Lambda \geq \frac{n-1}{n}$.
}

\vskip5mm
{\bf Proof.} Applying the inequality (1.3) to the matrix $A$ with only 
one non-zero entry on the $(i,j)$-place, we get
$$
\Var(X_i X_j) = \E X_i^2 X_j^2 - \delta_{ij} \leq \Lambda.
$$
Summing these bounds over all $i,j$ leads to $\E\,|X|^4 - n \leq n^2 \Lambda$.
But $\E\,|X|^4 \geq (\E\,|X|^2)^2 = n^2$. 
\qed

\vskip5mm
All the above definitions extend to complex-valued random variables $X_i$
using complex numbers $a_{ij}$ in the definition (1.3) (of course, $a_{ij}^2$ 
should be replaced with $|a_{ij}|^2$). Note that, if $\xi$ is 
a complex-valued random variable, its variance is defined by
$$
\Var(\xi) = \E\,|\xi - \E \xi|^2 = \E\, |\xi|^2 - |\E \xi|^2.
$$


\vskip5mm
\section{{\bf Classes of Distributions Satisfying
Second Order Correlation Condition}}
\setcounter{equation}{0}

Here we provide a few examples where functionals defined above may be
easily evaluated or properly estimated. Bounds are attained for the second
order correlation parameter for the following classes of distributions: 
i.i.d., coordinate-wise symmetric, log-concave and coordinate-wise symmetric,
and probability measures with a pectral gap.
 
As before, let $X = (X_1,\dots,X_n)$, 
$n \geq 2$. The case of independent components may be dealt with by simple 
calculation.

\vskip5mm
{\bf Proposition 3.1.} {\sl If the random variables 
$X_1,\dots,X_n$ are independent and have mean zero, then 
\be
\sigma_4^2(X) = \frac{1}{n}\,\sum_{i=1}^n \Var(X_i^2),
\en
\be
M_2(X) = \max_i \, \big(\E X_i^2\big)^{1/2}, \qquad
\Lambda(X) \leq 2 \max_i \,\E X_i^4.
\en
}

Note that equality (3.1) obviously extends to pairwise independent random 
variables with mean zero. The proof of the bound of $\Lambda(X)$ in (3.2) is 
similar to the one in Proposition 3.2 below, so we omit it.

Another class of illustrative examples is given by distributions of random 
vectors $X$ which are equal to $(\ep_1 X_1,\dots,\ep X_n)$ for arbitrary 
choices of signs $\ep_i = \pm 1$. We call such distributions coordinate-wise 
symmetric, although in the literature they are also called distributions with
unconditional basis. This class includes all symmetric product measures on 
$\R^n$ and corresponds to the case where the components $X_i$ are i.i.d. 
random variables with symmetric distributions on the line. It is therefore 
not surprising that many formulas like those in Proposition 3.1 extend to 
the coordinate-wise symmetric distributions. In particular, the first equality 
in (3.2) is still valid. As for $\Lambda(X)$, it
may be essentially reduced to the moment-type functional
$$
V(X) = 
\sup_{\theta \in S^{n-1}} \Var(\theta_1 X_1^2 + \dots + \theta_n X_n^2),
$$
representing the maximal eigenvalue of the matrix 
$\{{\rm cov}(X_i^2,X_j^2)\}_{i,j=1}^n$.

\vskip5mm
{\bf Proposition 3.2.} {\sl Given a random vector $X = (X_1,\dots,X_n)$ 
in $\R^n$ with a coor\-di\-nate-wise symmetric distribution, we have
\be
V(X) \leq \Lambda(X) \leq 2\max_i \, \E X_i^4 + V(X).
\en
If additionally the distribution of $X$
is invariant under permutations of coordinates, then
\be
\sigma_4^2(X) \, \leq \, \Lambda(X) \, \leq \, 2\, \E X_1^4 + \sigma_4^2(X),
\en
where the last term $\sigma_4^2(X)$ may be removed when
${\rm cov}(X_1^2,X_2^2) \leq 0$.
}

\vskip5mm
The proof of this proposition is rather elementary, but technical.
So, we postpone it to Section 9.

The following subfamily of coordinate-symmetric distributions admits 
a uniform bound on $\Lambda$. Let us recall that a (Borel) probability
measure $\mu$ on $\R^n$ is called log-concave, if it satisfies
the Brunn-Minkowski-type inequality
$$
\mu(tA + (1-t)B) \geq \mu(A)^t \mu(B)^{1-t}, \qquad 0 < t < 1,
$$
for all non-empty compact sets $A$ and $B$ in $\R^n$,
where $tA + (1-t)B \, = \, \{tx + (1-t)y: x \in A, \ y \in B\}$ denotes
the Minkowski weighted sum. An equivalent description
was given by Borell \cite{Bor}: the measure $\mu$ should be
supported on a closed convex set $V \subset \R$ and have a log-concave 
density $p$ with respect to the Lebesgue measure $\lambda_V$ on $V$ 
of the same dimension as $V$ (that is, $\log p$ is concave). Note that,
if $\mu$ is isotropic and log-concave, then necessarily $V$ has dimension $n$,
so that $\mu_V$ is the (full) Lebesgue measure.

\vskip5mm
{\bf Proposition 3.3.} {\sl Assume that the random vector $X$ in $\R^n$ is 
isotropic and has a coordinate-wise symmetric, log-concave distribution. 
Then
$$
\sigma_4^2(X) \leq \Lambda(X) \leq c.
$$ 
}

{\bf Proof.} The distribution of the random vector $(|X_1|,\dots,|X_n|)$ has
a log-concave, coordinate-wise non-increasing density. By a theorem due to 
Klartag \cite{K3}, the following weighted Poincar\'e-type inequality holds
$$
\Var(u(X)) \, \leq \, 4\,\E \sum_{i=1}^n X_i^2\,(\partial_i u(X_i))^2
$$
for any smooth even function $u$ on $\R^n$. Choosing
$u(x) = \theta_1 x_1^2 + \dots + \theta_n x_n^2$ with 
$\theta_1^2 + \dots + \theta_n^2 = 1$, we get
$$
\Var(u(X)) \, \leq \, 16 \sum_{i=1}^n \theta_i^2\, \E X_i^4 \, \leq \,
16\, \max_{i \leq n} \E X_i^4.
$$
In view of Proposition 3.2, we get
$$
\Lambda(X)  \, \leq \,
2\max_i \, \E X_i^4 + 
\sup_{\theta \in S^{n-1}} \Var(\theta_1 X_1^2 + \dots + \theta_n X_n^2) \\
 \, \leq \,
18\max_i \, \E X_i^4 .
$$

It remains to recall that $L^p$-norms of random variables with log-concave
distributions are equivalent to each other. In particular, for isotropic 
log-concave $X_i$'s, we have $\E X_i^4 \leq c\,(\E X_i^2)^2 = c$.
\qed

\vskip2mm
The above subclass may be potentially enlarged by considering the usual
Poincar\'e-type inequality
\be
\lambda_1 \Var(u(X)) \, \leq \, \E\, |\nabla u(X)|^2.
\en

\vskip2mm
{\bf Proposition 3.4.} {\sl Assume that a mean zero random vector $X$ in $\R^n$ 
satisfies a Poincar\'e-type inequality with constant $\lambda_1 > 0$. 
Then $M_2^2(X) \leq 1/\lambda_1$. Moreover,
$$
\sigma_4^2(X) \leq \Lambda(X) \leq \frac{4}{\lambda_1^2},
$$
and if $X$ isotropic, then
$$
\sigma_4^2(X) \leq \Lambda(X) \leq \frac{4}{\lambda_1}.
$$
}

{\bf Proof.}
Applying (3.5) to the linear functions $f(x) = \left<x,\theta\right>$,
$\theta \in S^{n-1}$, we obtain
$$
\lambda_1 \Var(\left<X,\theta\right>) \leq 1.
$$
If $X$ has mean zero, the latter means that 
$M_2^2(X) \leq 1/\lambda_1$. In particular, 
$\E\,X_j^2 \leq \frac{1}{\lambda_1}$. 
Taking the quadratic function 
$u(x) = \sum_{i,j=1}^n a_{ij}\, x_i x_j$ with $a_{ij} = a_{ji}$, 
we get, by Cauchy's inequality,
$$
\Var\bigg(\sum_{i,j=1}^n a_{ij} X_i X_j\bigg) \, \leq \,
\frac{4}{\lambda_1}\,\sum_{i=1}^n \E\, \bigg(\sum_{j=1}^n a_{ij} X_j\bigg)^2
 \, \leq \,
\frac{4}{\lambda_1}\,\sum_{i,j=1}^n a_{ij}^2\, \E X_j^2.
$$
Hence, the right-hand side does not exceed 
$4/\lambda_1^2$ subject to $\sum_{i,j=1}^n a_{ij}^2 \leq 1$,
and thus $\Lambda(X) \leq 4/\lambda_1^2$, while 
$\Lambda(X) \leq 4/\lambda_1$ in the isotropic case.
\qed


\vskip5mm
\section{{\bf Second Order Concentration on the Sphere}}
\setcounter{equation}{0}

Concentration of measure on the sphere means that the range of deviations of 
any Lipschitz function $u$ on the unit sphere $S^{n-1}$ is essentially of order 
at most $\frac{1}{\sqrt{n}}$, which may be strengthened as the subgaussian 
stochastic dominance $|u| \leq \frac{c}{\sqrt{n}}\,|Z|$ where $Z$ denotes 
a standard normal random variable (cf. \cite{M-S}, \cite{L}).
More precisely, there is a subgaussian deviation inequality
\be
\s_{n-1}\{|u(\theta)| \geq r\} \, \leq \, 2\, e^{-(n-1)\,r^2/2}, 
\qquad r > 0,
\en
valid whenever the smooth function $u$ has $\s_{n-1}$-mean zero and 
Lipschitz seminorm $\|u\|_{\rm Lip} \leq 1$. This may be partly seen
from the Poincar\'e inequality
\be
\int |u|^2\,d\s_{n-1} \leq \frac{1}{n-1}\int |\nabla u|^2\,d\s_{n-1}
\en
in the class of all smooth complex-valued $u$ with $\s_{n-1}$-mean zero. 
Although here there 
is equality for all linear functions, the spherical concentration phenomenon 
may be strengthened with respect to the dimension $n$ for a wide subclass 
of smooth functions. In order to facilitate applications, we shall not use
sphere intrinsic gradients but use Euclidean notions induced by the standard
embedding of the sphere. Here functions are defined in an open subset of $\R^n$ 
and their partial derivatives are understood in the usual sense.
We denote by $\nabla^2 u(x)$ the Hessian, that is, the $n \times n$ matrix
of second order partial derivative $\partial_{ij} u(x)$, and by
$I_n$ the identity $n \times n$ matrix. The next proposition summarizes
several recent results from \cite{B-C-G1}.

\vskip5mm
{\bf Proposition 4.1.} {\sl Suppose that a real-valued function $u$ is 
defined and $C^2$-smooth in some neighbourhood of $S^{n-1}$.
If $u$ is orthogonal to all affine functions in $L^2(\s_{n-1})$, then
\be
\int u^2\,d\s_{n-1} \leq \frac{5}{(n-1)^2}
\int \|\nabla^2 u - aI_n\|_{\rm HS}^2\,d\s_{n-1}
\en
for any $a \in \R$. Moreover, if $\|\nabla^2 u - a I_n\| \leq 1$ uniformly 
on $S^{n-1}$ for the operator norm, and the second integral 
in $(4.2)$ is bounded by $b$, then
\be
\int \exp\Big\{\frac{n-1}{2(1+4b)}\, |u|\Big\}\, d\s_{n-1} \leq 2.
\en
}

By Markov's inequality, (4.4) yields a corresponding large deviation 
bound, which may be stated informally as a subexponential stochastic dominance
$|u| \leq c_b\,(\frac{1}{\sqrt{n}} Z)^2$. In particular, this means that 
the deviations of $u$ are of order at most $1/n$.

The second order Poincar\'e-type inequality (4.3) obviously extends to all
complex-valued $u$ that are orthogonal to all affine functions on the sphere.
In this case, (4.4) may be applied separately to the real and imaginary
part of $u$, which results in
\be
\int \exp\Big\{\frac{n-1}{4(1+4b)}\, |u|\Big\}\, d\s_{n-1} \leq 2,
\en
assuming that $\|\nabla^2 u - a I_n\| \leq 1$ on $S^{n-1}$ 
for some $a \in {\mathbb C}$.


\vskip5mm
\section{{\bf Concentration of Characteristic Functions}}
\setcounter{equation}{0}

Given an isotropic random vector $X = (X_1,\dots,X_n)$ in $\R^n$,
introduce the following smooth functions
\be
u_t(\theta) = f_\theta(t) = \E\,e^{it\left<X,\theta\right>}, \qquad \theta \in \R^n,
\en
where $t \neq 0$ serves as a parameter. Note that,
for any fixed $\theta$, $t \rightarrow f_\theta(t)$ represents the 
characteristic function of the weighted sum $S_\theta = \left<X,\theta\right>$
with distribution function $F_\theta$, while the $\s_{n-1}$-mean of $u_t$,
$$
f(t) = \E_\theta f_\theta(t) = \E_\theta\, \E\,e^{it\left<X,\theta\right>},
$$
is the characteristic function of the avarage distribution function
\be
F(x) = \int F_\theta(x)\,d\s_{n-1}(\theta) = \E_\theta\, \P\{S_\theta \leq x\}, 
\qquad x \in \R. 
\en
Let us recall that we use $\E_\theta$ to denote integrals over the unit sphere
with respect to the uniform measure $\s_{n-1}$.

In order to study deviations of the functions $u_t$ from their $\s_{n-1}$-means
$f(t)$ on $S^{n-1}$, one may start from the Poincare inequality (4.2). Indeed,
differentiating the equality (5.1), we get that, for any $\theta' \in S^{n-1}$,
$$
\left<\nabla u_t(\theta),\theta'\right> = it \,
\E\left<X,\theta'\right> e^{it\left<X,\theta\right>},
$$
which, by Cauchy's inequality, implies
$$
|\left<\nabla u_t(\theta),\theta'\right>|^2 \leq t^2 \,
\E\left<X,\theta'\right>^2 = t^2.
$$
Taking the supremum over all $\theta'$, it follows that
$|\nabla u_t(\theta)| \leq |t|$, which means that $u_t$ has a Lipschitz
semi-norm $\|u_t\|_{\rm Lip} \leq |t|$ (on the whole space $\R^n$).
Therefore, by (4.2),
\be
\E_\theta\,|f_\theta(t) - f(t)|^2 \,\leq\, \frac{t^2}{n-1}.
\en

Thus, the deviations of $f_\theta(t)$ from $f(t)$ with respect to
$\theta \in S^{n-1}$ are of order at most $1/\sqrt{n}$ -- a property which 
may potentially be transfered to the analogous statement about the deviations 
of the distribution functions $F_\theta$ from $F$ in the sense of certain
weak metrics. 

In order to obtain better rates, we employ Proposition 4.1, assuming 
additionally that the random vector $X$ is symmetric and satisfies 
a second order correlatation condition (1.3) with parameter $\Lambda$.
To apply the bounds (4.3) and (4.5), we need to choose a suitable value 
$a \in \mathbb C$ and estimate the operator norm $\|\nabla^2 u_t - aI_n\|$ 
and the Hilbert-Schmidt norm $\|\nabla^2 u_t - aI_n\|_{\rm HS}$.
First note that, by further differentiation of (5.1), the Hessian of $u_t$ is given by
$$
\big[\nabla^2 u_t(\theta)\big]_{jk} =
\frac{\partial^2}{\partial \theta_j \partial \theta_k}\,f_\theta(t) =
-t^2\, \E\,X_j X_k\,e^{it\left<X,\theta\right>}
$$
for any fixed $t \in \R$.
Hence, a good choice could be $a = -t^2 f(t)$ in order to balance the diagonal
elements in the matrix of second derivatives of $u_t$. 
For any vector $v \in {\mathbb C}^n$ with complex components, using 
the canonical inner product in the complex $n$-space, we have
$$
\left<\nabla^2 u_t(\theta)v,v\right> = -t^2\,
\E\,|\left<X,v\right>|^2\,e^{it\left<X,\theta\right>}.
$$
Hence, with this choice of $a$, by the isotropy assumption,
$$
\big|\left<(\nabla^2 u_t(\theta) - aI_n)v,v\right>\big| \leq t^2\,
\E\,|\left<X,v\right>|^2 + |a|\,|v|^2 \leq 2t^2, \qquad |v|=1.
$$
This bound insures that 
\be
\|\nabla^2 u_t(\theta) - aI_n\| \leq 2t^2.
\en

In addition, putting $a(\theta) = -t^2 f_\theta(t)$, we have
\bee
\big\|\nabla^2 u_t(\theta) - a(\theta) I_n\big\|_{\rm HS}^2 
 & = &
\sum_{j,k = 1}^n \left|\nabla^2 u_t(\theta)_{jk} - a(\theta)\,\delta_{jk}\right|^2 \\
 & = &
\sup \bigg|\sum_{j,k=1}^n a_{jk}\,
\big(\nabla^2 u_t(\theta)_{jk} - a(\theta)\,\delta_{jk}\big)\bigg|^2 \\ 
 & = &
t^4 \sup \bigg|\,\E \sum_{j,k=1}^n a_{jk}\,\big(X_j X_k - \delta_{jk}\big)\,
e^{it\left<X,\theta\right>}\bigg|^2 \\
 & \leq &
t^4 \,\sup\,  \E\,\bigg|\sum_{j,k=1}^n a_{jk}\, (X_j X_k - \delta_{jk})\bigg|^2,
\ene
where the supremum is running over all complex numbers $a_{jk}$ such that 
$\sum_{j,k=1}^n |a_{jk}|^2 = 1$. But, under this constraint (with complex 
coefficients), due to the second order correlation condition, the last
expectation is bounded by $\Lambda$, so that
\be
\big\|\nabla^2 u_t(\theta) - a(\theta) I_n\big\|_{\rm HS}^2 \, \leq \,
\Lambda t^4
\en
for all $\theta$. On the other hand, by (5.3),
\be
\E_\theta\, \big\|(a(\theta) - a) I_n\big\|_{\rm HS}^2 \, = \, 
n t^4\,\E_\theta\,|f_\theta(t) - f(t)|^2 \, \leq \, 2t^6,
\en
since $\frac{n}{n-1} \leq 2$. The two bounds give
\be
\E_\theta\, \big\|\nabla^2 u_t(\theta) - aI_n\big\|_{\rm HS}^2 \, \leq \, 
2\Lambda t^4 + 4t^6.
\en

Note that (5.6) is worse in comparison with (5.5) in the variable $t$. Nevertheless,
applying the second order Poincar\'e-type inequality, it is possible to
improve the resulting inequality (5.7) for reasonably long $t$-intervals.
Since the distribution of $X$ is symmetric about the origin, the characteristic 
functions $f_\theta(t)$ are even with respect to $\theta$, i.e., 
$f_{-\theta}(t) = f_\theta(t)$. Hence, they are orthogonal in the Hilbert space 
$L^2(\s^{n-1})$ to all linear functions on the sphere. Thus, the conditions 
of Proposition 4.1 are fulfilled 
for the function $u = u_t$, and using (5.7), the inequality (4.3) gives
\be
\E_\theta\,|f_\theta(t) - f(t)|^2 \, \leq \, \frac{5}{(n-1)^2}\,(2\Lambda t^4 + 4t^6).
\en
This bound allows us to improve (5.6) to the form
\bee
\E_\theta\, \big\|(a(\theta) - a) I_n\big\|_{\rm HS}^2 
 & \leq & 
nt^4\,\frac{5}{(n-1)^2}\,(2\Lambda t^4 + 4t^6) \\ 
 & \leq & 
\frac{40\,\Lambda t^4}{n}\,(t^4 + 2t^6) \, \leq \, 120\,\Lambda t^4,
\ene
where in the last inequality we assume that $|t| \leq n^{1/6}$. Combining
this with (5.5), we therefore obtain that
$$
\E_\theta\, \big\|\nabla^2 u_t(\theta) - aI_n\big\|_{\rm HS}^2 \, \leq \, 
c\Lambda t^4.
$$
In view of (4.3), this already gives the inequality (5.9) below.

To get a stronger deviation inequality, let us recall (5.4), so that to 
conclude that the conditions of Proposition 4.1 (in its second part) are 
fulfilled for the function
$$
u(\theta) = \frac{1}{2t^2}\,(f_\theta(t) - f(t)), \qquad \theta \in \R^n, \ \ 
0 < |t| \leq n^{1/6},
$$
with  parameter $b = c\Lambda$ (which bounded away from zero).
Applying (4.5), we arrive at:

\vskip5mm
{\bf Corollary 5.1.} {\sl Let $X$ be an isotropic random vector in $\R^n$ with 
a symmetric distribution and finite constant $\Lambda$. Then the characteristic 
functions~$f_\theta(t) = \E\,e^{it\left<X,\theta\right>}$ satisfy
\be
\E_\theta\, |f_\theta(t) - f(t)|^2\,\leq\,\frac{c}{n^2} \Lambda t^4
\en
whenever $0 < |t| \leq n^{1/6}$. Moreover,
\be
\E_\theta\, \exp\Big\{\frac{n}{c \Lambda t^2}\, |f_\theta(t) - f(t)|\Big\} 
\leq 2.
\en
}

As we have seen, removing the constraint $|t| \leq n^{1/6}$, (5.9) may be
replaced with a weaker inequality (5.8). When applying the latter to the estimation
of $\rho(F_\theta,F)$ via Lemma 6.1 below, we would gain an additional
$\log n$ factor in Theorem 1.1.

\vskip5mm
\section{{\bf Proof of Theorem 1.1}}
\setcounter{equation}{0}

Based on the deviation inequalities (5.9)-(5.10), Fourier analytic tools yield 
bounds for the closeness of the distribution functions $F_\theta$ to the 
$s_{n-1}$-mean distribution function $F$ defined in (5.2). The following 
Berry-Esseen-type bound can be found in \cite{B-C-G3}, cf. Lemma 6.2, 
which we state in the case $p=2$.

\vskip5mm
{\bf Lemma 6.1.} {\sl Suppose that a random vector $X$ in $\R^n$ has 
a finite moment of order $4$, with $\E\,|X|^2 = n$. Then, for all 
$T \geq T_0 > 0$,
\begin{eqnarray}
c\,\E_\theta\, \rho(F_\theta,F)
 & \leq & 
\int_0^{T_0} \E_\theta\,|f_\theta(t) - f(t)|\,\frac{dt}{t} \nonumber \\
 & & + \ 
\frac{M_4^4 + \sigma_4^2}{n}\,\bigg(1 + \log\frac{T}{T_0}\bigg)
+ \frac{1}{T} + e^{-T_0^2/16}.
\end{eqnarray}
}

{\bf Proof of Theorem 1.1}.
Applying Propositions 2.1-2.2 and using the isotropy assumption, we have 
$M_4^4 + \sigma_4^2 \leq 1 + 2\Lambda \leq 4\Lambda$. Hence, (6.1) yields
\bee
c_1\,\E_\theta\, \rho(F_\theta,F) 
 & \leq &
\int_0^{T_0} \E_\theta\,|f_\theta(t) - f(t)|\,\frac{dt}{t} \\
 & & + \ 
\frac{\Lambda}{n}\,\bigg(1 + \log\frac{T}{T_0}\bigg)
+ \frac{1}{T} + e^{-T_0^2/16}.
\ene
Here, the integrand may be estimated by virtue of (5.9), and then we get
$$
c_2\,\E_\theta\, \rho(F_\theta,F) \, \leq \, \frac{1}{n}\, T_0^2 \sqrt{\Lambda} +
\frac{1}{n}\,\bigg(1 + \log\frac{T}{T_0}\bigg) \Lambda
+ \frac{1}{T} + e^{-T_0^2/16},
$$
provided that $T_0 \leq n^{1/6}$.
As a natural choice, take $T_0 = 5\sqrt{\log n}$, $T = 5n$
(assuming that $n$ is large enough), which leads to the bound
\be
c_3\,\E_\theta\, \rho(F_\theta,F) \, \leq \, \frac{\log n}{n}\,\Lambda.
\en
We finally refer to \cite{B-C-G2}, 
Theorem 1.1, cf. also \cite{B-C-G3}, Corollary 4.2, where the estimate
\be
\rho(F,\Phi) \leq c\,\frac{1 + \sigma_4^2}{n}
\en
was derived. Using $\sigma_4^2 \leq \Lambda$ and combining (6.2) 
with the triangle inequality for $\rho$,
we arrive at the desired inequality (1.4). 
\qed

\vskip5mm
{\bf Remark 6.2.} Under proper moment assumptions and using the spherical 
deviation inequality (5.10), one may derive large deviation bounds for 
$\rho(F_\theta,\Phi)$ as well. In particular, suppose that
\be
\E\,e^{|S_\theta|/\beta} \leq 2 \qquad 
\en
for all $\theta \in S^{n-1}$ with some $\beta>0$. Then, in the setting
of Theorem 1.1,
$$
\s_{n-1}\Big\{\rho(F_\theta,F) \geq 
\frac{c\log n}{n}\,(\Lambda + \beta^4)\, r\Big\} \, \leq \, 
cn\,\exp\{-r^{1/8}\}, \qquad r \geq 0.
$$
In other words, with high $\s_{n-1}$-probability, 
$$
\rho(F_\theta,F) \leq \frac{c\,(\log n)^9}{n}\,(\Lambda + \beta^4).
$$
For details we refer the interested reader to \cite{B-C-G4}.

\vskip5mm
\section{{\bf The log-concave case}}
\setcounter{equation}{0}

Specializing to the class of isotropic log-concave distributions on $\R^n$,
first let us comment on the unconditional statement with a standard rate 
of normal approximation as indicated in the inequality (1.7). If the isotropic 
random vector $X$ has a uniform distribution over a symmetric convex body 
in $\R^n$, it was shown by Anttila, Ball, and Perissinaki that
\be
\s_{n-1}\big\{\rho(F_\theta,F) \geq r\big\} \, \leq \, 
c_1 \sqrt{n}\,\log n\,e^{-c_2 nr^2}, \qquad r > 0
\en
(actually with $c_2 = 50$, cf. \cite{A-B-P}). With a different argument, this inequality 
has been extended to arbitrary isotropic log-concave distributions in \cite{B1}.
In both papers, as a main step, it was observed that, for every point $x \in \R$, 
the function $u(\theta) = F_\theta(x)$ has a bounded Lipschitz semi-norm
on the unit sphere, so that one may apply the spherical
concentration inequality (4.1), leading to
$$
\s_{n-1}\big\{|F_\theta(x) - F(x)| \geq r\big\} \, \leq \, 
2\,e^{-c nr^2}, \qquad r > 0.
$$
Since $\rho(F_\theta,F) \leq 1$, (7.1) readily yields an upper bound
$$
\E_\theta\, \rho(F_\theta,F) \leq c\sqrt{\frac{\log n}{n}}.
$$
Combining it with (6.3) and applying the triangle inequality for the metric $\rho$,
we therefore obtain the normal approximation on average in the form of the relation
\be
c\, \E_\theta\, \rho(F_\theta,\Phi) \, \leq \, \sqrt{\frac{\log n}{n}} +
\frac{\sigma_4^2}{n}.
\en
It remains to involve the bound $\sigma_4^2 \leq c\sqrt{n}$, which
was recently derived by Lee and Vempala \cite{L-V}, and then we arrive at (1.7).

A thin-shell conjecture, raised in \cite{B-K},
asserts that the functional $\sigma_4^2(X)$, or equivalently $\Var(|X|)$,
is actually bounded by a dimension-free (and thus universal) constant over the whole
class of isotropic log-concave random vectors $X$ in $\R^n$. Specializing
to the convex body case, a similar concentration hypothesis was also
suggested in \cite{A-B-P}. It states that the deviation inequality
$$
\P\Big\{\Big|\frac{|X|}{\sqrt{n}} - 1\Big| \geq \ep_n\Big\} \leq \ep_n
$$
holds true with $\ep_n \leq c\,(\log n)/\sqrt{n}$.  The boundedness of
$\sigma_4^2$ allows one to take a slightly thinner shell with $\ep_n = c/\sqrt{n}$.
Anyhow, the bound (7.2) subject to the thin-shell conjecture still leads to the
standard normal approximation as in (1.7).

Note that, by the Poincar\'e-type inequality (1.5) applied with $u(x) = |x|^2$, 
one gets $\sigma_4^2 \leq 4/\lambda_1$, so that the thin-shell conjecture
is formally weaker than the K-L-S (which is further precised in Proposition 3.4). 
On the other hand, recently Eldan \cite{E}
has developed a new localization technique, in essense reducing
the stronger hypothesis to the weaker one modulo a logarithmic factor.
It is is therefore possible to state Corollary 1.2 alternatively as follows.

\vskip5mm
{\bf Corollary 7.1.} {\sl Let $X$ be an isotropic random vector in $\R^n$
with a symmetric log-concave distribution. Assuming that 
the thin-shell conjecture is true, we have
$$
\E_\theta\, \rho(F_\theta,\Phi) \leq \frac{c\,(\log n)^3}{n}.
$$
}

{\bf Proof.} Combining Theorem 1.1 with Proposition 3.4, we get
\be
\E_\theta\, \rho(F_\theta,\Phi) \, \leq \, \frac{c}{\lambda_{1,n} n} \log n,
\en
where $\lambda_1 = \lambda_{1,n}$ is the smallest spectral gap in 
the Poincar\'e-type inequality over the class of all isotropic log-concave 
probability measures on $\R^n$. Assuming the K-L-S conjecture, 
$\lambda_{1,n}$ is bounded away from zero, which thus leads to the inequality (1.6)
of Corollary 1.2.
Within the same class, this quantity may be related to the largest value 
$\sigma_{4,n}^2 = \sup_X \sigma_4^2(X)$. Namely, as shown by Eldan \cite{E}, 
\be
\frac{1}{\lambda_{1,n}} \, \leq \, 
c \log n \sum_{k=1}^n \frac{\sigma_{4,k}^2}{k}.
\en
In particular, the bound of the form
$\sigma_{4,n}^2 \leq c_1 n^\alpha$ $(0 \leq \alpha \leq 1)$ implies that
\be
\lambda_{1,n}^{-1} \leq c\, \eta_\alpha(n)
\en
with $\eta_\alpha(n) = \frac{c_1}{\alpha} \,n^\alpha\,\log n$ in case
$\alpha>0$ and $\eta_0(n) = 3c_1\,(\log n)^2$. 
It remains to apply (7.5) in (7.3) with $\alpha = 0$.
\qed


\section{{\bf From the normal approximation to the shin shell}}
\setcounter{equation}{0}

\vskip5mm
To refine the relationship between the central limit theorem and
the thin-shell problem, let us complement Corollary 7.1 by the following
general statement involving the maximal $\psi_1$-norm of linear
functionals of $X$. 

\vskip5mm
{\bf Proposition 8.1.} {\sl Let $X$ be a random vector in $\R^n$
with $\E\,|X|^2 = n$, satisfying the moment condition $(6.4)$
with some $\beta>0$. Then
\be
c\,\sigma_4^2(X) \, \leq \,
n\,(\beta \log n)^4\ \E_\theta\, \rho(F_\theta,\Phi) + 
\frac{\beta^4}{n^4} + 1.
\en
}

In the isotropic log-concave case, the condition (6.4) is fulfilled
with some absolute constant $\beta$ (by the well-known Borell's 
Lemma 3.1 in \cite{Bor}), and this simplifies (8.1) to
$$
c\,\sigma_4^2(X) \, \leq \,
n\,(\log n)^4\ \E_\theta\, \rho(F_\theta,\Phi) + 1.
$$
Hence, the potential property 
\be
\E_\theta\, \rho(F_\theta,\Phi) \leq \frac{c\,\log n}{n}
\en
as in Corollary 1.2 would imply that
\be
\sigma_4^2(X) \leq c\, (\log n)^5,
\en 
assuming additionally that the distribution of $X$ is symmetric about zero. 
But, the symmetry condition may easily be dropped. Indeed, define 
$X' = (X-Y)/\sqrt{2}$, where $Y$ is an independent copy of a random 
vector $X$ with an isotropic log-concave distribution on $\R^n$.
Then, the distribution of $X'$ is isotropic, log-concave, and symmetric 
about zero. Moreover,
\bee
\sigma_4^2(X') 
 & = &
\frac{1}{2n}\,\Var\big(|X|^2 + |Y|^2 - 2\left<X,Y\right>\big) \\ 
 & = & 
\frac{1}{2n}\,\Var(|X|^2) + \frac{1}{2n}\,\Var(|Y|^2) + 
\frac{2}{n}\, \E \left<X,Y\right>^2 \, = \, 
\sigma_4^2(X) + 2.
\ene
Hence, once (8.3) is true for the random vector $X'$, it continues to hold 
for $X$ as well (with other constant).

Note also that, applying Eldan's inequality (7.4) together with (8.3), 
from the normal approximation  (8.2) we get
$$
\lambda_{1,n}^{-1} \leq c\,(\log n)^7.
$$

{\bf Proof of Proposition 8.1.} In view of the triangle inequality 
$\rho(F,\Phi) \leq \E_\theta\, \rho(F_\theta,\Phi)$,
it is sufficient to derive (8.1) for $\rho(F,\Phi)$ in place of 
$\E_\theta\, \rho(F_\theta,\Phi)$. This means that we need in essense 
to reverse the inequality (6.3) by using (6.4). To this aim, let us rewrite 
the definition (5.2) as
$$
F(x) = \P\{|X|\, \theta_1 \leq x\}, \qquad x \in \R,
$$
where we assume that $X$ and $\theta = (\theta_1,\dots,\theta_n) \in S^{n-1}$ 
(as a random vector uniformly distributed on the sphere)
are independent. This description yields
$$
\int_{-\infty}^\infty x^4\,dF(x) \, = \,
\E\, |X|^4\ \E_\theta\, \theta_1^4 \, = \,
(n^2 + \sigma_4^2 n)\,\frac{3}{n(n+2)},
$$
or equivalently
\be
\int_{-\infty}^\infty x^4\,dF(x) - \int_{-\infty}^\infty x^4\,d\Phi(x)
 = \frac{3}{n+2}\,(\sigma_4^2 - 2),
\en
where $\sigma_4^2 = \sigma_4^2(X)$.
On the other hand, it follows from (6.4) that
$$
\int_{-\infty}^\infty e^{|x|/\beta}\,dF(x) \leq 2.
$$
Using $t^2 \leq 4e^{-2}\,e^t$ ($t \geq 0$) together with the property
$\int_{-\infty}^\infty x^2\,dF(x) = \frac{1}{n}\,\E\,|X|^2 = 1$, we
have $\beta \geq e/\sqrt{8}$, which can be used to derive the bounds 
$$
1 - \Phi(x) \leq \frac{1}{2}\,e^{-x^2/2} \leq 2\,e^{-x/\beta}, \qquad x \geq 0.
$$
In addition, by Markov's inequality,
$F(-x) + (1 - F(x)) \leq 2\,e^{-x/\beta}$, so that
$$
|F(-x) - \Phi(-x)| + |F(x) - \Phi(x)| \, \leq \, 6\,e^{-x/\beta}
$$
for all $x \geq 0$. Hence, integrating by parts, 
we see that, for any $T \geq 6\beta$, the left-hand side
of (8.4) does not exceed in absolute value
\bee
4 \int_{-T}^T |x|^3\,|F(x) - \Phi(x)|\,dx + 24
\int_T^\infty x^3\,e^{-x/\beta}\,dx & & \\
 & & \hskip-40mm \leq \
2T^4\,\rho(F,\Phi) + 48\,\beta T^3\,e^{-T/\beta}.
\ene
Choosing $T = 9\beta \log n$ and recalling (8.4), we get 
$$
\sigma_4^2 \, \leq \ 6 + cn\, \Big[
\beta^4\,(\log n)^4\,\rho(F,\Phi) + \frac{\beta^4}{n^9}\,(\log n)^3\Big].
$$
\qed


\section{{\bf Proof of Proposition 3.2}}
\setcounter{equation}{0}

The lower bound on $\Lambda$ in (3.3) immediately follows 
from (1.3) by choosing the coefficients to be of the form 
$a_{ij} = \theta_i \delta_{ij}$. For the upper bound, put 
$v_i^2 = \E X_i^2$ and define
$$
X^{(2)}_{ij} = X_i X_j - \E X_i X_j = X_i X_j - \delta_{ij} v_i^2.
$$
The covariances of these mean zero random variables are given by
\begin{eqnarray}
\E\,X^{(2)}_{ij} X^{(2)}_{kl}
 & = & 
\E\,(X_i X_j - \delta_{ij} v_i^2) \, X_k X_l \\
 & = &
\E\,X_i X_j X_k X_l - \delta_{ij} \delta_{kl}\, v_i^2 v_k^2.  \nonumber
\end{eqnarray}

Case 1: $i \neq j$. By the symmetry with respect to the coordinate axes, 
the right-hand side of (9.1) is vanishing unless 
$(i,j) = (k,l)$ or $(i,j) = (l,k)$. In both cases, it is equal to
$$
\E\,X^{(2)}_{ij} X^{(2)}_{ij} = \E\,X^{(2)}_{ij} X^{(2)}_{ji} =
\E X_i^2 X_j^2.
$$

Case 2: $i = j$. The right-hand side in (9.1) is non-zero only 
when~$k = l$. 

Case 2a): $i = j$, $k = l$, $i \neq k$.  The right-hand side in (9.1) 
is equal to
$$
\E\,X^{(2)}_{ii} X^{(2)}_{kk} = \ \E X_i^2 X_k^2 - \E X_i^2 \, \E X_k^2 = 
{\rm cov}(X_i^2,X_k^2).
$$

Case 2b): $i = j = k = l$. The right-hand side is equal to
$$
\E\,X^{(2)}_{ii} X^{(2)}_{ii} = \E X_i^4 - \E X_i^2 \, \E X_i^2 = \Var(X_i^2).
$$
In both subcases, $\E\,X^{(2)}_{ii} X^{(2)}_{kk} = {\rm cov}(X_i^2,X_k^2)$.
Therefore, for any collection of real numbers $a_{ij}$ such that $a_{ij} = a_{ji}$
and $\sum_{i,j=1}^n a_{ij}^2 = 1$,
\bee
\Var\bigg(\sum_{i,j=1}^n a_{ij} X_i X_j\bigg) 
 & = &
\sum_{i,j = 1}^n \sum_{k,l = 1}^n a_{ij}\, a_{kl}\, \E\,X^{(2)}_{ij} X^{(2)}_{kl} \\
 & = &
2\ \sum_{i \neq j} a_{ij}^2 \, \E X_i^2 X_j^2 + 
\sum_{i,k} a_{ii}\, a_{kk}\, {\rm cov}(X_i^2,X_k^2).
\ene
Here, the first sum on the right-hand side does not exceed 
$$
\max_{i \neq j}\, \E X_i^2 X_j^2\, \sum_{i \neq j} a_{ij}^2 \, \leq \, 
\max_{i}\, \E X_i^4\, \sum_{i \neq j} a_{ij}^2 \, \leq \, 
\max_{i}\, \E X_i^4
$$
(by applying Cauchy's inequality). As for the second sum, 
it does not exceed $V(X)$, and we obtain 
\be
\Lambda(X) \, \leq \, 2\,\max_{i\neq j} \, \E X_i^2 X_j^2 + V(X),
\en
from which (3.3) follows immediately.

As for (3.4), recall that the first inequality always holds, 
cf. Proposition 2.1. For the second one, let us note that
\be
\sigma_4^2(X) \, = \, \frac{1}{n}\,\Var(|X|^2) \, = \, 
\Var(X_1^2) + (n-1)\, {\rm cov}(X_1^2,X_2^2)
\en
and that, for any 
$\theta = (\theta_1,\dots,\theta_n) \in S^{n-1}$,
$$
\Var(\theta_1 X_1^2 + \dots + \theta_n X_n^2) \, = \,
\Big(\sum_{i =1}^n \theta_i\Big)^2\, {\rm cov}(X_1^2,X_2^2) -
{\rm cov}(X_1^2,X_2^2) + \Var(X_1^2).
$$
Here, in the case ${\rm cov}(X_1^2,X_2^2) \geq 0$, the right-hand side is 
maximized for equal coefficients, and recalling (9.3), we then get
$$
\Var(\theta_1 X_1^2 + \dots + \theta_n X_n^2) \leq 
(n-1){\rm cov}(X_1^2,X_2^2) + \Var(X_1^2) = \sigma_4^2(X).
$$
Hence, (9.2) implies (3.4).
In the case ${\rm cov}(X_1^2,X_2^2) \leq 0$, we similarly conclude that
\bee
\Var(\theta_1 X_1^2 + \dots + \theta_n X_n^2) 
 & \leq & 
-{\rm cov}(X_1^2,X_2^2) + \Var(X_1^2) \\
 & = & 
\E X_1^4 - \E X_1^2 X_2^2,
\ene
which means that $V(X) \leq \E X_1^4 - \E X_1^2 X_2^2$. Thus, by (9.2),
$$
\Lambda(X) \, \leq \, 2\, \E X_1^2 X_2^2 + V(X) \leq \E X_1^2 X_2^2 + \E X_1^4
 \leq 2\,\E X_1^4.
$$
Hence, (3.4) follows in this case as well even without the
$\sigma_4^2(X)$-functional.
\qed


\vskip5mm
\section{{\bf Historical Remarks}}
\setcounter{equation}{0}

Finally, let us give a short overview on results related to Theorem 1.1
(some account can also be found in the book \cite{B-G-V-V}).
It is natural to distinguish between two types of results.

\vskip5mm
{\bf 10.1. Deviations of $F_\theta$ from the mean distribution $F$ in different 
metrics.} The paper by Sudakov [Su] starts with the hypothesis
$$
\E\,\Big(\sum_{i=1}^n a_i X_i\Big)^2 \, \leq \, M_2^2\, \sum_{i=1}^n a_i^2,
\qquad a_i \in \R,
$$
which may be called a first order correlation condition. Here, an optimal value 
$M_2 = M_2(X)$ is the same functional we considered in 
Section 2; equivalently, $M_2^2$ represents the maximal eigenvalue 
of the correlation operator for the random vector $X$. As was shown in \cite{Su},
if $M_2$ is bounded, and $n$ is large, then most of $F_\theta$ are close to the
average distribution $F$ in the sense of the Kantorovich or $L^1$-distance
$$
W_1(F_\theta,F) = \|F_\theta - F\|_{L^1(\R,dx)} = \int_{-\infty}^\infty
|F_\theta(x) - F(x)|\,dx.
$$
A closely related observation was also made by Diaconis and Freedman \cite{D-F}.
A somewhat different scheme, in which the coefficient vectors are 
drawn from the Gaussian measure $\mu_n$ on $\R^n$ 
with mean zero and covariance matrix $\frac{1}{n}\,I_n$, was also
considered by Nagaev \cite{N} and von Weizs\"acker \cite{W}. In particular, 
assuming that $M_1 = 1$, [N] contains a quantitative bound 
\be
\bigg(\int \|F_\theta - F\|_{L^2(\R,dx)}^2\,d\mu_n(\theta)\bigg)^{1/2} 
\leq \frac{1}{(\pi n)^{1/4}}
\en
for the $L^2$-distance between the distribution functions. When 
the coefficients have a special structure, similar phenomena were
considered in \cite{B2}, \cite{B-G}.

Returning to the spherical measure $\s_{n-1}$, the rate as in (10.1)
is achieved for the L\'evy distance as well. More precisely, there is 
a general bound
$$
\E_\theta\, L(F_\theta,F) \leq c \,\frac{\log n}{n^{1/4}},
$$
where the constant $c$ depends on $M_1$ only, cf. \cite{B4}. Large deviation 
bounds on $L(F_\theta,F)$ were given in \cite{B1} in the isotropic case. 
As was already discussed in Section 7,
the rate and deviation bounds may be essentially improved and be 
stated for the stronger Kolmogorov distance, when the random vector 
$X$ has an isotropic log-concave distribution.

Quantitative variants of Sudakov's theorem for $W_1$ were studied in \cite{B3}, 
where it was shown that, for any $p>1$,
$$
\E_\theta\, W_1(F_\theta,F) \, \leq \, \frac{12\,p}{p-1}\,M_p\,n^{-\frac{p-1}{2p}}.
$$
The rate is thus approaching $1/\sqrt{n}$ for growing $p$.
Under a stronger assumption (6.4), the above inequality easily implies
$$
\E_\theta\, W_1(F_\theta,F) \leq c\beta \,\frac{\log n}{\sqrt{n}}.
$$
Here, the logarithmic term may be removed, if $X$ has an isotropic 
log-concave distribution (by virtue of Proposition 3.1 in \cite{B1}).
Note that, in all these results, the rates are not better than a multiple of
$1/\sqrt{n}$.

\vskip5mm
{\bf 10.2. Deviations of $F_\theta$ from the standard normal distribution 
function $\Phi$.} 
To study the approximation of $F_\theta$ by the standard normal 
distribution function for most of $\theta$'s, one is led to determine
rates for the distance $\rho(F,\Phi)$, which may be reduced to the estimation
of $\sigma_4^2(X)$ (via relation (6.3)). In fact, the control of the 
two functionals, $M_4$ and $\sigma_4$, is sufficient to guarantee a standard 
rate of normal approximation for $F_\theta$ on average. 
As was shown in \cite{B-C-G3}, we have
$$
\E_\theta\, \rho(F_\theta,\Phi) \, \leq \, 
c\,(M_4^3 + \sigma_4^{3/2})\,\frac{1}{\sqrt{n}}.
$$

Note that Theorem 1.1 essentially improves this bound as long as $\Lambda$
is of the same order as $M_4$ and $\sigma_4$. However, whether or not
$\Lambda$ and even $\sigma_4$ is bounded might be a difficult problem for 
some classes of distributions on $\R^n$ such as the class of
isotropic log-concave probability measures. For this class, the property that
$\rho(F_\theta,\Phi)$ is small for most of $\theta$ (when $n$ is large)
was first established by Klartag \cite{K1}. In particular,
$\E_\theta\, \rho(F_\theta,\Phi) \leq \ep_n \rightarrow 0$
as $n \rightarrow \infty$ uniformly over the class.
For further refinements in this direction, see \cite{K2}, \cite{E-K}, \cite{G-M}.

There is also a number of results, where the coefficients are fixed, and
$\rho(F_\theta,\Phi)$ are bounded by a quantity, which depends on $\theta$ 
as well, cf. e.g. \cite{M}, \cite{M-M}. One striking result due to Klartag \cite{K3}
should be mentioned: If  the random vector $X$ in $\R^n$ is isotropic 
and has a coordinate-wise symmetric, isotropic log-concave distribution, then 
$$
\rho(F_\theta,\Phi) \, \leq \, c \sum_{k=1}^n \theta_k^4.
$$ 
Moreover, a similar bound holds true for the stronger total variation distance.
This is of course more precise in comparison
with the average estimate $\E_\theta \rho(F_\theta,\Phi) \leq c/n$.

\section*{Acknowledgment}
\hskip -3mm The authors would like to thank the two referees for careful
reading of the manuscript and valuable comments.

\end{document}